\documentclass[12pt]{article}

\usepackage{amssymb}
\usepackage{amsmath}
\usepackage{amsbsy}
\usepackage{amscd}
\usepackage{amsfonts}
\usepackage{amsthm}
\usepackage{mathrsfs}
\usepackage{verbatim}
\usepackage[colorlinks]{hyperref}
\usepackage{fullpage}
\usepackage{mathdots}
\usepackage{graphicx,subfigure}
\usepackage[english]{babel}
\usepackage[utf8]{inputenc}
\usepackage{tikz}
\usepackage{adjustbox}
\usepackage{authblk}
\usepackage{caption}

\usepackage{mathtext}

\usepackage{tikz-cd}
\usetikzlibrary{backgrounds,fit, matrix}
\usetikzlibrary{positioning}
\usetikzlibrary{calc,through,chains}
\usetikzlibrary{arrows,shapes,snakes,automata, petri}
\usepackage[boxsize=1.25em, centerboxes]{ytableau}

\newtheorem{Theorem}[subsection]{Theorem}
\newtheorem{Corollary}[subsection]{Corollary}
\newtheorem{Lemma}[subsection]{Lemma}
\newtheorem{Proposition}[subsection]{Proposition}

\theoremstyle{definition}
\newtheorem{Definition}[subsection]{Definition}
\newtheorem{Example}[subsection]{Example}

\newtheorem{Notation}[subsection]{Notation}

\newtheorem{Remark}[subsection]{Remark}

\numberwithin{equation}{section}

\newcommand{\Z}{{\mathbb Z}}

\newcommand{\N}{{\mathbb N}}

\begin{document}

\title{Incidence Monoids: Automorphisms and Complexity}

\author[1]{Mahir Bilen Can}

\affil[1]{\normalsize{
Tulane University, New Orleans, Louisiana\\
mahirbilencan@gmail.com}} 

\maketitle

\begin{abstract}
The algebraic monoid structure of an incidence algebra is investigated. 
We show that the multiplicative structure alone determines the algebra automorphisms of the incidence algebra. 
We present a formula that expresses the complexity of the incidence monoid with respect to the two sided action of its maximal torus in terms of the zeta polynomial of the poset. In addition, we characterize the finite (connected) posets whose incidence monoids have complexity $\leq 1$. 
Finally, we determine the covering relations of the adherence order on the incidence monoid of a star poset. 
\\

\\


\end{abstract}

\section{Introduction}

Let $k$ be an algebraically closed field. 
If $R$ is a finite dimensional $k$-algebra, 
to what extent the group of algebra automorphisms of $R$ is controlled by its underlying multiplicative semigroup structure? We answer this question for a sufficiently self-interesting class of $k$-algebras, 
namely, the incidence algebras of finite posets. 

Let $P$ be a finite poset. The {\em incidence algebra of $P$ over $k$} is the associative $k$-algebra defined by 
$I(P):= \{ f: P\times P \to k :\ f(x,y) = 0 \text{ if $x \nleq y$}\}$.
The monoid structure on $I(P)$ is given by the convolution product, 
$(f\cdot g)(x,y) = \sum_{x \leq z \leq y} f(x,z) g(z,y)$ for $f,g\in I(P)$ and $x,y,z\in P$.
The group of invertible elements of this monoid is solvable. 
In the sequel, when we want to emphasize the monoid structure of $I(P)$, we will refer to it as the {\em incidence monoid of $P$}.
The importance of incidence algebras for combinatorics stems from the fact that for two (locally finite) posets $P$ and $Q$, $I(P)$ and $I(Q)$ are isomorphic as $k$-algebras if and only if $P$ and $Q$ are isomorphic as posets~\cite{Stanley1970}. 

Let $M$ be an irreducible linear algebraic monoid with zero. 
As a semigroup, $M$ is regular if and only if the Zariski closure of the solvable radical of the unit group $G$ of $M$ is a completely regular semigroup.
Such regular monoids are reductive, and conversely, every reductive monoid is regular~\cite{Putcha,Renner}. 
The normal reductive monoids~\cite{Renner1985,Vinberg1995,Rittatore1998} as well as normal irreducible completely regular algebraic monoids~\cite{Renner1989} are classified according to their toric data. 
Every associative algebra over a field has an algebraic monoid structure~\cite{Okninski2014},
but, in general, such an algebraic monoid may not be regular in the semigroup sense. 
Likewise, in general, the incidence algebras are neither regular nor completely regular. 
Nevertheless, as we will show in the sequel, in every incidence algebra, once the unipotent radical of the group of units is fixed, there is a spectrum of irreducible, smooth, and completely regular algebraic submonoids. 
In fact, the extreme elements of this spectrum are, ordered by inclusion, form a lower meet semilattice. 
For some posets, we are able to determine the spectrum explicitly.

We go back to our first question. 
On one hand, the automorphism groups of incidence algebras are well-studied objects~\cite{Stanley1970,Backlawski1972,Spiegel2001,DrozdKolesnik}.
For an excellent exposition of automorphism groups of incidence algebras, we recommend the article~\cite{BrusamarelloLewis}.
On the other hand, not much is known about the automorphism group of linear algebraic monoids $M$ except for the regular monoids with zero, see~\cite{Putcha1983}. One of our goals in this paper is to contribute to the study of the automorphism groups of general algebraic monoids. 
In this regard, our first main theorem states that the algebraic monoid structure alone determines the $k$-algebra automorphisms of an incidence algebra.

\begin{Theorem}\label{T:main1}
Let $I(P)$ be the incidence monoid of a finite poset $P$.
Then the group of algebraic semigroup automorphisms, denoted by $\text{Aut}_{s.g.}(I(P))$, is equal to 
the group of $k$-algebra automorphisms of $I(P)$.
\end{Theorem}

Our second goal is to establish some further connections between combinatorics and the theory of equivariant group embeddings. 
Let $G$ be a connected linear algebraic group, and let $B=UT$ be a Borel subgroup of $G$. 
Here, $U$ is the unipotent radical of $B$, and $T$ is the maximal torus contained in $B$. 
If $X$ is an irreducible variety on which $G$ acts morphically, then let us agree to indicate it by the notation $G:X$. 
One of the most important invariants of $G:X$ is called the {\em complexity}, denoted by $c_G(X)$.
It is defined as the transcendence degree of $k(X)^B$ over $k$, where $k(X)$ is the field of rational functions on $X$ and $k(X)^B$ is the $B$-invariants. 
For a detailed study of this concept, we recommend the monographs by Grosshans~\cite{Grosshans} and Timashev~\cite{Timashev}.

We now turn back to the incidence monoids. 
Let $P$ be a finite poset together with a linear extension $P= \{x_1,\dots, x_n\}$. 
As we mentioned before, the $I(P)$ is a solvable monoid, that is, the unit group $G(P)$ of $I(P)$ is a solvable algebraic group. 
The doubled unit group, $G(P)\times G(P)$ acts on $I(P)$ by $(g_1,g_2)\cdot x = g_1 x g_2^{-1}$ for every $g_1,g_2 \in G(P)$ and $x\in I(P)$. 
Clearly, the orbit of the identity element of $I(P)$ is equal to $G(P)$, hence, we have an open dense orbit of $G(P)\times G(P)$. 
It follows that the complexity of $G(P)\times G(P): I(P)$ is always 0, so, in this case there is not much ado. 
Actually, this is true even if we view $I(P)$ as a $G(P)$-variety via the left, or the right multiplication action of $G(P)$.
Let $T$ be a maximal torus of $G(P)$. We restrict the action of $G(P)\times G(P)$ to the subgroup $T\times T$. 
Now the question of computing the complexity of $T\times T: I(P)$ becomes a combinatorially interesting problem.

The {\em zeta polynomial} of $P$ is the unique polynomial function on $\Z$, denoted by $Z(P,n)$, such that if $n\geq 2$, 
then $Z(P,n)$ is the number of multichains $t_1 \leq t_2 \leq \cdots \leq t_{n-1}$ in $P$.

\begin{Theorem}\label{T:main2}
Let $P$ be a finite connected poset. Then we have 
\begin{align*}
c_{T\times T}(I(P)) = Z(P,3)-2Z(P,2)+1.
\end{align*}
\end{Theorem}

By using Theorem~\ref{T:main2}, we are able to classify all finite (connected) posets whose incidence monoid has complexity at most 1.

\begin{Theorem}\label{T:main3}
Let $P$ be a finite connected poset.
Then the following statements hold: 
\begin{enumerate}
\item The complexity $c_{T\times T}(I(P))$ is zero if and only if 
all maximal chains of $P$ has length 1, and the underlying graph of Hasse diagram of $P$ has no circuits. 
\item The complexity $c_{T\times T}(I(P))$ is one if and only if $P$ is a graded poset, and one of the following statements holds:
\begin{enumerate}
\item all maximal chains of $P$ have length 1, and the underlying graph of the Hasse diagram of $P$ has exactly one circuit; 
\item $P$ has a unique maximal chain of length 2, and the underlying graph of the Hasse diagram of $P$ has no circuits. 
\end{enumerate}
\end{enumerate}
\end{Theorem}

Let $S$ be a torus. An irreducible normal $S$-variety with complexity zero is a toric variety. 
Hence, part 1. of our Theorem~\ref{T:main3} is a classification theorem for toric incidence monoids $I(P)$, where $P$ is a connected poset.
We note that, by combining Theorems~\ref{T:main1} and~\ref{T:main3},
we obtain a concrete description of the semigroup automorphism group of an incidence monoids of complexity $\leq 1$.

Let $R$ denote the set of orbits for an action $H:X$, where $H$ is a connected algebraic group. 
The order on $R$ that is determined by the inclusions of orbit closures is called the {\em adherence order}.
If $H= B\times B$ and $X$ is a reductive monoid, where $B$ is a Borel subgroup of $X$, 
then the adherence order is called the {\em Bruhat-Chevalley-Renner} order. 
In this case, $R$ has the structure of a finite monoid~\cite{Renner1985} 
since it is given by the quotient $\overline{N_G(T)}/T$, where $N_G(T)$ is the normalizer of the maximal torus $T$ in
the unit group $G$ of $X$. 
This finite monoid is called the {\em Renner monoid} of $X$ and it contains the Weyl group $N_G(T)/T$ as its unit group. 
For the special case $X:=\mathbf{Mat}_n$, that is the monoid of $n\times n$ matrices, 
the Renner monoid is called the {\em rook monoid} since it consists of 
0/1 matrices of order $n$ such that in every row and column, there exists at most one 1.
For a concrete description of the Bruhat-Chevalley-Renner order on the Renner monoid of $\mathbf{Mat}_n$, see~\cite{CanRenner,Can2019}. 
Although for a broad class of spherical varieties (including reductive monoids) there is a general description of the adherence order for the Borel subgroup actions~\cite{PPR,Timashev1994}, for $B$-varieties which are not necessarily spherical, the adherence order of Borel orbits remains largely unexplored. 
In our next result, we present some results in that direction for the action $G(P)\times G(P):I(P)$, where $I(P)$ is a certain toric incidence monoid.

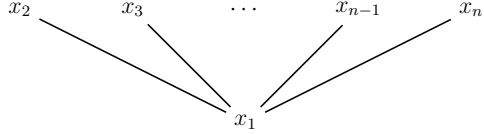
\begin{figure}[htp]
\begin{center}
\scalebox{.75}{
\begin{tikzpicture}[scale=1]

\node (1) at (0,-1) {$x_1$};
\node (2) at (-4,1) {$x_2$};
\node (3) at (-2,1) {$x_3$};
\node (4) at (0,1) {$\cdots$};
\node (5) at (2,1) {$x_{n-1}$};
\node (6) at (4,1) {$x_n$};
\draw[thick,-] (1) to (2); 
\draw[thick,-] (1) to (3);
\draw[thick,-] (1) to (5); 
\draw[thick,-] (1) to (6); 
 
\end{tikzpicture}
}
\end{center}

\caption{The star poset on $n$ vertices.}
\label{F:toric1}
\end{figure}

Let $P_n$ denote the {\em star poset}, which is defined by its Hasse diagram as shown in Figure~\ref{F:toric1}.

\begin{Theorem}\label{T:main4}
Let $P_n$ denote the star poset, where $n \geq 2$. Let $G(P_n)$ denote the unit group of $I(P_n)$,
and let $T$ denote the maximal torus of $G(P_n)$. Then we have the following results:
\begin{enumerate}
\item The incidence monoid of $P_n$ is a smooth toric variety of dimension $2n-1$.
\item The group of algebraic semigroup automorphisms of $I(P_n)$ is isomorphic to $G(P_n)\rtimes S_{n-1}$, where $S_{n-1}$ is the $(n-1)$-th symmetric group.
\item The adherence order for $T\times T : I(P)$ is isomorphic to the Boolean algebra of 
subsets of an $(2n-1)$-element set. 
In particular, there are $2^{2n-1}$ $T\times T$-orbits in $I(P_n)$. 
\item We have concrete descriptions of the rank function as well as the covering relations of the adherence order for 
$G(P_n)\times G(P_n) : I(P_n)$.
In particular, there are $2^{n-1}+3^{n-1}$ $G(P_n)\times G(P_n)$-orbits in $I(P_n)$. 
\end{enumerate}
\end{Theorem}

Finally, we want to mention that some of our results have interesting generalizations in the context of incidence monoids of finite preordered sets (quasi-ordered sets).

The organization of our paper is as follows. 
In Section~\ref{S:Preliminaries}, we introduce basic poset terminology that we will use in the rest of our paper.
In Subsection~\ref{SS:LAM}, we review some basic facts about linear algebraic groups and monoids. 
Actually, we prove our first result (Proposition~\ref{T:Clifford}) in this subsection, which is about 
completely regular algebraic submonoids of an incidence monoid. 
We show that for every antichain (clutter) of $P$, 
there is an irreducible, smooth, and completely regular algebraic submonoid whose unipotent radical is equal to the unipotent radical of $I(P)$.  
Furthermore, we observe that the set of all such completely regular submonoids of $I(P)$ is closed under intersections. 
In Section~\ref{S:Automorphisms}, we prove our first main result Theorem~\ref{T:main1}.
As a consequence of this theorem we see that the outer automorphism group of an incidence monoid $I(P)$ is 
determined by the topology of the Hasse diagram $\mathscr{H}_P$ and the automorphism group of the poset $P$. 

In Section~\ref{S:Complexity} we prove Theorems~\ref{T:main2} and~\ref{T:main3}.
In addition, we compute the automorphism groups of incidence monoids of complexity $\leq 1$.
The purpose of Section~\ref{S:Star} is to analyze the adherence order on star posets. 
We prove Theorem~\ref{T:main4} in that section. 
Furthermore, we present some applications of our previous results to the incidence monoid $I(P_n)$. 


\section{Preliminaries}\label{S:Preliminaries}

We begin with setting up basic poset terminology. 

\subsection{Posets.}\label{SS:Posets}

In this article we are concerned with finite posets only. 
Therefore, we will not mention it again.

A {\em chain} in $P$ is a totally ordered subposet of $P$. 
An {\em antichain} in $P$ is a subset $Q\subset P$ such that no two elements of $Q$ are comparable. 
Let $x$ and $y$ be two elements from $P$. 
If $x< y$ and $x\leq z < y$ implies that $z=x$, then $y$ is said to {\em cover} $x$.
In this case, $(x,y)$ is called a covering relation.  
A sequence of elements $C=(x_1,\dots, x_n)$ such that $x_1 \leq x_2 \leq \cdots \leq x_n$ is called a {\em chain}. 
If, in addition, we have $x_i \neq x_{i+1}$ for every $i\in \{1,\dots, n-1\}$, then $C$ is called a {\em saturated chain}.
A {\em maximal chain} is a saturated chain that is not a subsequence of any chain other than itself. 
The {\em length} $\ell(C)$ of a saturated chain $C=(x_1,\dots, x_n)$ is defined as $n-1$. 
A poset $P$ is said to be {\em bounded} if there exists a positive integer $m$ such that 
for every saturated chain $C$ of $P$, we have $\ell(C) < m$.

A poset $P$ is called a {\em graded (or ranked) poset} if every maximal chain in $P$ has the same length. 
In this case, up to a translation by an integer, there is a unique function $\rho: P \to \N$ such that $\rho(y) = \rho(x) +1$ whenever $y$ covers $x$.
We will refer to $\rho$ as the {\em rank function} of $P$. 
If $P$ has a unique minimal (resp. maximal) element, then we will denote it by $\hat{0}$ (resp. by $\hat{1}$).

The {\em Hasse diagram of $P$}, denoted by $\mathscr{H}_P$, is the directed graph whose vertices are the elements of $P$. 
In this graph, there is a directed edge between $x$ and $y$, directed towards $y$, whenever $y$ covers $x$. 
The information that is provided by $\mathscr{H}_P$ uniquely determines $P$. 
A poset $P$ is said to be {\em connected} if its Hasse diagram is connected. 
Clearly, if a poset possesses a maximum or a minimum element, then it is connected. 
A subposet $Q$ of $P$ is called a {\em connected component} if $\mathscr{H}_Q$ is a connected component of $\mathscr{H}_P$. 
Note that a Hasse diagram is a {\em simple graph} in the sense that there is at most one edge between any two vertices. 
A {\em walk} in a graph is a finite or infinite sequence of edges which joins a sequence of vertices.
A {\em trail} is a walk in which all edges are distinct.
A {\em circuit} is a non-empty trail in which the first vertex is equal to the last vertex.

We now briefly mention some basic facts about the incidence algebras. 
First of all, the incidence algebra of a poset $P$ can be defined over an arbitrary unital associative algebra, however, since we are interested in their geometry, 
we focus on the incidence algebras that are defined over algebraically closed fields.

The origins of the theory of incidence algebras go back to Dedekind,
but it was Rota~\cite{Rota1964} who first used them for a systematic study of the M\"obius inversion formula. 
Let us briefly review this development in our way to setup some notation that we will use later. 

Let $\mathbf{B}_n$ denote the standard Borel subgroup of all invertible upper triangular $n\times n$ matrices with entries from $k$. 
Then its Zariski closure, $\overline{\mathbf{B}_n}$ in $\mathbf{Mat}_n$ has the structure of not only an algebraic monoid but also a finite dimensional $k$-algebra. 
A {\em linear extension of $P$} is a labeling of the elements of $P$ by $\{1,\dots, n\}$, where $|P|=n$, in such a way that 
if $x_i \leq x_j$, then $i \leq  j$. 
A linear extension of $P$ enables us to define a natural $k$-algebra embedding,
\begin{align}\label{A:linearextension}
\varphi: I(P) &\longrightarrow \overline{\mathbf{B}_n} \notag \\
f &\longmapsto (f_{ij})_{i,j=1}^n,
\end{align}
where, for every $i,j\in \{1,\dots, n\}$, we have  
\begin{align*}
f_{ij} :=
\begin{cases}
f(x_i,x_j) & \ \text{ if $x_i \leq x_j$,}\\
0 & \ \text{ otherwise}.
\end{cases}
\end{align*}
Under this embedding, the {\em zeta function of $P$}, which is defined by  
\[
\zeta (x,y) := 
\begin{cases}
1 &\ \text{ if $x\leq y$,}\\
0 & \ \text{ otherwise,}
\end{cases}
\] 
is mapped to a unipotent matrix in $\overline{\mathbf{B}_n}$.
Therefore, it is a unit in $I(P)$. 
The inverse of $\zeta$ is called the {\em M\"obius function of $P$}, denoted by $\mu$. 
The M\"obius function can be computed recursively as follows: 
\begin{align*}
\mu(s,s) &= 1,\ \text{  for all $s\in P$},\\
\mu(s,u) &= -\sum_{s\leq t < u } \mu(s,t),\ \text{ for all $s<u$ in $P$}.
\end{align*}

\begin{Example}
Let $P$ denote the poset whose Hasse diagram is as in Figure~\ref{F:firstexample}. 
\begin{figure}[htp]
\begin{center}
\scalebox{.75}{
\begin{tikzpicture}[scale=1.5]

\node (0) at (0,0) {$x_1$};
\node (1) at (-2,1) {$x_2$};
\node (2) at (2,1) {$x_3$};
\node (3) at (-2,2) {$x_4$};
\node (4) at (2,2) {$x_5$};

\draw[thick,-] (0) to (1);
\draw[thick,-] (0) to (2);
\draw[thick,-] (1) to (3);
\draw[thick,-] (1) to (4);
\draw[thick,-] (2) to (3);
\draw[thick,-] (2) to (4);

\end{tikzpicture}
}

\end{center}
\caption{The Hasse diagram of a poset.}
\label{F:firstexample}
\end{figure}
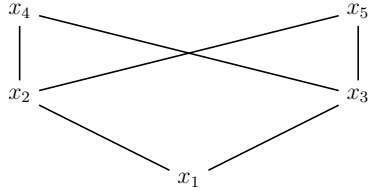
Notice that the subscripts of the vertices of $\mathscr{H}_P$ give a linear extension of $P$. 
Then the images of the zeta and the M\"obius functions of $P$ under the embedding 
$\varphi : I(P) \to \overline{\mathbf{B}_5}$ are given by 
\begin{align*}
\varphi(\zeta) = 
\begin{bmatrix}
1 & 1 & 1 & 1 & 1 \\
0 & 1 & 0 & 1 & 1 \\
0 & 0 & 1 & 1 & 1 \\
0 & 0 & 0 & 1 & 0 \\
0 & 0 & 0 & 0 & 1 \\
\end{bmatrix}
\qquad \text{and} \qquad
\varphi(\mu)= 
\begin{bmatrix}
1 & -1 & -1 & 1 & 1 \\
0 & 1 & 0 & -1 & -1 \\
0 & 0 & 1 & -1 & -1 \\
0 & 0 & 0 & 1 & 0 \\
0 & 0 & 0 & 0 & 1 \\
\end{bmatrix}.
\end{align*}
\end{Example}


\subsection{Linear Algebraic Monoids.}\label{SS:LAM}

We begin with basic algebraic group theory terminology. 
Let $G$ be a linear algebraic group. 
The maximal normal connected unipotent subgroup of $G$, denoted by $\mathscr{R}_u(G)$,
is called the {\em unipotent radical of $G$}. 
If $\mathscr{R}_u(G)$ is the trivial subgroup, then $G$ is said to be {\em reductive}.
The {\em radical of $G$} is the maximal normal connected solvable subgroup; it is denoted by $\mathscr{R}(G)$.
Clearly, $\mathscr{R}_u(G)$ is a subgroup of $\mathscr{R}(G)$. 
A Borel subgroup of $G$ is a maximal connected solvable subgroup. 
Any two Borel subgroups in $G$ are conjugate. 
If $G$ is connected, then its radical is the intersection of all Borel subgroups of $G$. 
The {\em standard Borel subgroup} of the $n$-th general linear group $\mathbf{GL}_n$ is denoted by $\mathbf{B}_n$. 
The maximal diagonal torus in $\mathbf{B}_n$ will be denoted by $\mathbf{T}_n$, and the unipotent radical of $\mathbf{B}_n$
will be denoted by $\mathbf{U}_n$. 
The group of characters of a group $H$ will be denoted by $\mathcal{X}(H)$.

If $I(P)$ is the incidence monoid of a poset $P$, then we will denote its unit group by $G(P)$. 
Then we have a decomposition, 
\[
G(P) = T \ltimes U_P,
\]
where $T$ is isomorphic to $\mathbf{T}_n$, and $U_P$ is the unipotent radical, $U_P:=\mathscr{R}_u(G(P))$.
Note that, the coordinates of $U_P$ are completely determined by $P$,
\[
U_P = \{ (u_{ij})_{i,j=1}^n \in \mathbf{U}_n :\ \text{$u_{ij}= 0$ if $x_i \nleq x_j$} \}.
\]

Let $M$ be an irreducible linear algebraic monoid defined over an algebraically closed field $k$. 
Let $G$ denote the unit group of $M$, and let $E(M)$ denote the variety of idempotents of $M$. 
If $G$ is a reductive group, then $M$ is said to be a {\em reductive monoid}. 
In this case, we have $M = G E(M) = E(M)G$. 
A semigroup $M$ is called {\em regular} if for every $a\in M$ there exists $x\in M$ such that $axa = a$.  
As we mentioned before, by a deep result of Putcha and Renner, we know that an irreducible linear algebraic monoid $M$ with zero
is regular if and only if $G$ is reductive; for a proof, see any of the monographs~\cite{Putcha,Renner}. 
If $M$ does not have a zero, then the reductiveness of $M$ is equivalent to the following two conditions~\cite{Huang2001}:
1) $M$ is regular, 2) the semigroup kernel of $M$ is reductive. 
However, note that, every incidence monoid has a zero.

A {\em completely regular semigroup} (or, a {\em Clifford semigroup}) is a semigroup in which every element lies in a subgroup. 
By a result of Putcha, we know that an irreducible algebraic semigroup $S$ with zero is completely regular if and only if the unit group of $S$ is a torus~\cite[Theorem 5.12]{Putcha}. 
Therefore, any incidence monoid of a poset $P$ which has a connected component $Q\subseteq P$ 
with $|Q|\geq 2$ cannot be completely regular. 

Next, we will discuss some constructions for building completely regular monoids from an incidence monoid. 
First of all, in every incidence monoid $I(P)$, the closure of a maximal torus $T$ in $I(P)$ gives a completely regular submonoid, whose semisimple rank is $|P|$.  But these completely regular submonoids do not capture any essential information about $P$ other than its cardinality. 
Let $n$ denote $|P|$ so that $\dim T = n$. 
Let $\varepsilon_1,\dots, \varepsilon_n$ be the coordinate functions on the algebraic monoid $\overline{T}$. 
For $A\subseteq \{1,\dots, n\}$, we define the associated {\em antichain torus}, denoted by $S_A$, as  
\[
S_A = \{ x \in T:\ \varepsilon_i(x) = 1\ \text{ for every $i\notin A$}\}.
\]
The Zariski closure of $S_A$ in $\overline{T}$ is a commutative monoid.
Its character monoid given by 
\[
\mathcal{X}(\overline{S_A}) = \bigoplus_{i\in  A} \N \varepsilon_i. 
\]

Let $e$ be the minimal idempotent from $E(\overline{S_A})$. 
Then $e$ is the zero element of $\overline{S_A}$. 
Following~\cite{Renner1989}, we define three unipotent subgroups of $U_P$:
\begin{align*}
U_+ &:= \{ u\in U_P:\ eu =e \}, \\
U_- &:= \{ u\in U_P:\ ue =e \}, \\
U_0 &:= \{ u\in U_P:\ eu =eu =e \}.
\end{align*}

We are now ready to present our construction which associates to an antichain a completely regular algebraic monoid. 

\begin{Proposition}\label{T:Clifford}
Let $P$ be a poset. Let $\{x_1,\dots, x_n\}$ be a linear extension of $P$. 
Let $A$ be a subset of $\{1,\dots, n\}$ such that $P_A:= \{x_i :i\in A\}$ is an antichain of $P$. 
Then the algebraic submonoid defined by 
\[
I(P)_A := \overline {S_A \ltimes U_P} \subseteq  I(P)
\]
is an irreducible, smooth, completely regular submonoid of $I(P)$. 
Furthermore, if $A$ is a maximal antichain in $P$, then $I(P)_A$ is maximal with respect to these properties. 
\end{Proposition}

\begin{proof}
As before, we denote by $T$ the maximal diagonal torus $\mathbf{T}_n$ in $\mathbf{B}_n$. 
Let $t$ be an element of $S_A$ such that $t= \text{diag}(t_1,\dots, t_n) \in \mathbf{T}_n$. 
Then $t_i =1$ for every $i\notin A$. 
Since $U_P$ is a subgroup of $\mathbf{U}_n$, its Lie algebra is a Lie subalgebra of $\text{Lie}(\mathbf{U}_n)$
which spanned as a $k$-vector space by the elementary matrices $E_{i,j}$, $1\leq i < j \leq n$. 
To compute the weights of $S_A$ on $\text{Lie}(U_P)$, let $i,j\in \{1,\dots, n\}$ be two indices such that 
the $(i,j)$-th elementary matrix $E_{i,j}$ is a basis element for $\text{Lie}(U_P)$.
Then $x_i \leq x_j$ in $P$, and therefore, we have one of the following three cases:
1) $\{i,j\}\cap A= \{i\}$, 2) $\{i,j\}\cap  A= \{j\}$, 3) $\{i,j\} \cap A =\emptyset$. 
In the first case (resp. the second case), the action of $t$ on $E_{i,j}$ is given by $t_i E_{i,j}$ (resp. by $t_j^{-1} E_{i,j}$). 
In the third case, it is given by $E_{i,j}$.
Therefore, the corresponding weights of the $S_A$-action are given by $\varepsilon_i$, $-\varepsilon_j$, or $0$ in the order of the three cases that we listed above.

Now we will use a result of Renner~\cite[Converse of Theorem 2.6]{Renner1989}:
Suppose that we are given the following data:
\begin{enumerate}
\item A torus action $S:U$, where $U$ is a unipotent algebraic group. 
\item A normal torus embedding $S\to \overline{S}$ such that $0\in \overline{S}$ and 
the weights of $S:U$ are contained in $\mathcal{X}(\overline{S})\cup -\mathcal{X}(\overline{S})$.
\end{enumerate}
Then there exists a unique structure of an irreducible normal completely regular algebraic monoid on 
$U_+ \times \overline{S} \times C_U(S) \times U_-$ extending the group law on $S \ltimes U$.  
Here, $C_U(S)$ denotes the centralizer of $S$ in $U$. 

We just checked that our subtorus $S_A$ and the unipotent subgroup $U_P$ together with $\overline{S_A}$ in $\overline{T}$
satisfy these conditions, hence, by Renner's theorem, $I(P)_A $ is an irreducible normal completely regular algebraic monoid.
The maximality follows from the fact that $A$ is a maximal subset that ensures that 
the weights of the conjugation action of $S_A$ on $U_P$ are all contained in $\mathcal{X}(\overline{S_A})\cup -\mathcal{X}(\overline{S_A})$.
This finishes the proof of our assertion. 

\end{proof}

We want to emphasize the fact that $I(P)_A$ does not have a zero;
the zero of $\overline{S_A}$ is not a zero for $I(P)_A$. 
Notice also that the unipotent radical $U_P$ is a connected closed subset of $I(P)$, so, it is automatically an irreducible, smooth, completely regular algebraic submonoid. 
The completely regular monoids of the form $I(P)_A$ as in Proposition~\ref{T:Clifford} behave well under the intersections.
The following corollary is now a consequence of the definitions and Proposition~\ref{T:Clifford}. 
\begin{Corollary}\label{C:semi}
Let $P$ be as in Proposition~\ref{T:Clifford}.
If $A$ and $B$ are antichains in $P$, then we have 
\begin{enumerate}
\item $A \subseteq B \iff \text{$I(P)_A$ is a submonoid of $I(P)_B$}$,
\item $I(P)_{A\cap B} = I(P)_A \cap I(P)_B$.
\end{enumerate}
In particular, the poset $(\{I(P)_A :\ \text{$A$ is an antichain in $P$}\}, \subseteq)$ is a finite meet-semilattice with minimum 
element $U_P$. 
\end{Corollary}
The meet-semilattice that we defined in Corollary~\ref{C:semi} will be called the {\em intersection lattice of antichains}.

We go back to our review of the algebraic monoids. 
We already mentioned in the introduction that $\mathbf{B}_n\times \mathbf{B}_n$-orbits in $\mathbf{Mat}_n$ are parametrized by 
the rook monoid, denoted by $R_n$, which consists of 0/1 matrices $x\in \mathbf{Mat}_n$ such that 
there exists at most one 1 in each row and column of $x$. 
We will denote by $id_n$ (resp. by $0$), the $n\times n$ identity matrix (resp. the zero of $\mathbf{Mat}_n$).

A {\em partial transformation on $\{1,\dots, n\}$} is just a function $g: A\to \{1,\dots, n\}$, where $A\subseteq \{1,\dots, n\}$. 
A partial transformation $g: A\to B$ is said to be an {\em injective partial transformation} if $g$ is injective.
Let $Q_n$ denote the set of all partial transformations on $\{1,\dots, n\}$. 
There is a useful faithful representation of $Q_n$ in $\mathbf{Mat}_n$ which is defined as follows: 
For $g\in Q_n$, the corresponding matrix $\rho(g)$ in $\mathbf{Mat}_n$ has a 1 in the $(i,j)$-th entry if $g(i) = j$;
otherwise it has a 0.
Note that, the image of the set of injective partial transformations under $\rho : Q_n\to \mathbf{Mat}_n$ 
is precisely the rook monoid $R_n$. 
In the sequel, it will also be useful to represent the elements of $Q_n$ in {\em word notation}: 
A word $a_1\dots a_n$ with entries from $\{0,1,\dots, n\}$ represents $g\in Q_n$ if $a_i = j$ whenever $g(i)=j$,
and $a_i=0$ if $g(i)$ is undefined.

The rook monoid has a natural partial order, called the {\em Bruhat-Chevalley-Renner order}, 
\[
x \leq y \iff \mathbf{B}_n x \mathbf{B}_n \subseteq \overline{\mathbf{B}_n y \mathbf{B}_n}\ \qquad (x,y\in R_n).
\]
Then the lower interval $[0,id_n]$ in $(R_n,\leq)$ consists of rook matrices that represents the 
$\mathbf{B}_n\times \mathbf{B}_n$-orbits in $\overline{\mathbf{B}_n}$.
For a detailed study of the interval $[0,id_n]$, see~\cite{CC}. 
It turns out that some elements of $[0,id_n]$ are relevant to our work.  

\begin{Lemma}\label{L:someofthem} 
Let $P$ be a poset with a linear extension $P=\{x_1,\dots, x_n\}$. 
Let $x$ be a rook matrix from the interval $[0,id_n]$ in $R_n$.
Then $G(P)x G(P)$ is contained $I(P)$ if and only if $x\in I(P)$. 
Furthermore, if $x$ and $y$ are two rooks from $I(P)\cap R_n$, then $x\neq y$ if and only if $G(P)x G(P) \cap G(P)y G(P) = \emptyset$.
\end{Lemma}

\begin{proof}
The first if and only if statement does not need a proof. 
For the second if and only if statement, it suffices to observe the following facts: 
1) $G(P) \subseteq \mathbf{B}_n$, 2) $[0,id_n]$ consists of $\mathbf{B}_n\times \mathbf{B}_n$-orbit representatives in $\overline{\mathbf{B}_n}$, 3) two orbits are either equal or disjoint. 
\end{proof}

\section{Automorphisms}\label{S:Automorphisms}

In this section, we will prove our first main Theorem~\ref{T:main1}.
This theorem will provide us with a concrete description of the outer automorphisms of the incidence monoids. 
As we mentioned before, a similar result is obtained by Putcha for regular monoids~\cite{Putcha1983}. 
In the next section, we will discuss toric algebraic monoids.
Then we will use Corollary~\ref{C:maincorollary} of Theorem~\ref{T:main1} to compute their automorphism groups.

Let $K$ be a finite dimensional $k$-algebra. 
As an algebraic variety, $K$ is isomorphic to the affine space $\mathbb{A}_k^{\dim K}$. 
In particular, $K$ has the structure of a smooth, connected linear algebraic monoid. 
We denote by $\text{Aut}_{alg.}(K)$ (resp. by $\text{Aut}_{s.g.}(K)$) the group of $k$-algebra automorphisms of $K$ (resp. the group of algebraic semigroup automorphisms of $K$).
If $P$ is a finite poset, we will denote by $\text{Aut}(P)$ the group of automorphisms of $P$.

\begin{proof}[Proof of Theorem~\ref{T:main1}]

Let $P$ be a poset with a linear extension $P=\{x_1,\dots, x_n\}$. 
Let $K$ denote its incidence algebra, $K:=I(P)$, and let $G$ denote the unit group of $K$. 
We will show that $\text{Aut}_{alg.}(K) \cong \text{Aut}_{s.g.}(K)$.
Clearly, any $k$-algebra automorphism is also a semigroup automorphism, so, we have 
$\text{Aut}_{alg.}(K) \subseteq \text{Aut}_{s.g.}(K)$.
Conversely, let $\varphi$ be an algebraic semigroup automorphism of the monoid $K$. 
Without loss of generality we will assume that $G$ is of the form $G=T\ltimes U$,
where $T$ is the maximal diagonal torus $\mathbf{T}_n$, and $U$ is a subgroup of $\mathbf{U}_n$. 
Since $\varphi$ restricts to give an automorphism of $G$, we know that $\varphi(T)$ is a maximal torus of $G$. 
By composing $\varphi$ with an inner automorphism, we may assume that it maps $T$ isomorphically onto $T$. 
Since it is an algebraic semigroup automorphism, $\varphi$ maps the closure of $T$ isomorphically onto the closure of $T$, 
$\varphi(\overline{T}) = \overline{T}$. 
Let $E(\overline{T})$ denote the set of idempotents of $\overline{T}$.
Let $e$ be an element from $ E(\overline{T})$. 
Then the equalities $\varphi(e) = \varphi(e^2) = \varphi(e) \varphi(e)$ imply that $\varphi$ defines an automorphism of $E(\overline{T})$. 
We let $E_1$ denote the set of rank 1 idempotents in $E(\overline{T})$. 
We label the elements of $E_1$ in such a way that they match the linear extension $\{x_1,\dots, x_n\}$ on $P$:
\[
E_1:=\{e_{x_1},\dots, e_{x_n}\}.
\]
Note that the partial order on $P$ is equivalent to the following partial order on $E_1$: 
\begin{align}\label{A:partialorderalternative}
e_{x_i} \leq e_{x_j} \iff e_{x_i} K e_{x_j} \neq 0.
\end{align}
Since $\varphi$ is an automorphism, we see that 
\[
e_{x_i} K e_{x_j} \neq 0 \iff 
\varphi(e_{x_i} K e_{x_j}) \neq 0 \iff
\varphi (e_{x_i}) K \varphi(e_{x_j}) \neq 0.
\]
In other words, $\varphi$ preserves the partial order on $P$.
Hence, $\varphi \in \text{Aut}(P)$. 

What we have shown so far did not take into consideration the following subgroup of $\text{Aut}_{s.g.}(K)$:
\[
 \text{Aut}_{s.g.}(K,\overline{T} ) := \{ \psi \in \text{Aut}_{s.g.}(K):\ \psi ( y) =y\ \text{ for all } y\in \overline{T} \}.
\]
Note that $K$ is spanned, as a $k$-algebra, by the elementary matrices $E_{i,j}$, where $i,j\in \{1,\dots, n\}$ are such that $e_{x_i} \leq e_{x_j}$.
Note also that $\text{span}_k \{ E_{i,j} \} =   e_{x_i} K e_{x_j}$ for every $i,j\in \{1,\dots, n\}$.
Therefore, if $\psi \in  \text{Aut}_{s.g.}(K,\overline{T} )$, then 
\[
\psi (E_{i,j}) = c(i,j) E_{i,j}\ \text{ for some $c(i,j) \in k^*$}.
\]
But these automorphisms are precisely the $k$-algebra automorphisms that are called the {\em multiplicative automorphisms}, see~\cite[Section 7.3]{SO}. 
Hence, we showed that every algebraic semigroup automorphism of $K$ is a $k$-algebra automorphism.
This finishes the proof of our theorem.
\end{proof}

For a linear algebraic monoid $M$ with unit group $G$,
to emphasize the semigroup structure, let us denote by $\text{Inn}_{s.g.}(M)$ the group of inner automorphisms of $M$, 
\[
 \text{Inn}_{s.g.}(M) := \{ \varphi \in \text{Aut}_{s.g.}(M):\ \text{there exists $g\in G$ s.t. $\varphi(h) = ghg^{-1}$ for every $h\in M$}\}.
\]
In a similar way, let us denote by $\text{Out}_{s.g.}(M)$ the group of outer algebraic semigroup automorphisms,  
\[
\text{Out}_{s.g.}(M) := \text{Aut}_{s.g.}(M)/\text{Inn}_{s.g.}(M).
\]

\begin{Corollary}\label{C:maincorollary}
If $K$ is the incidence monoid of a poset $P$, then its group of outer semigroup automorphisms is given by 
\begin{align}\label{A:cohom}
 \text{Out}_{s.g.}(K) = H^1( \mathscr{H}_P, k^\times) \rtimes \text{Aut}(P),
\end{align}
where $\text{Aut}(P)$ is the automorphism group of $P$, $H^1$ is the first cohomology group of $\mathscr{H}_P$ viewed as a simplicial complex. 
\end{Corollary}

\begin{proof}
By~\cite[Theorem 4]{DrozdKolesnik} we know that the outer automorphism group of an incidence algebra of a finite poset $P$ defined over $k$ 
is exactly the right hand side of (\ref{A:cohom}).
Clearly, the group of inner algebraic semigroup automorphisms of $K$ is equal to the group of inner automorphisms of $K$ as a $k$-algebra.
Therefore, by Theorem~\ref{T:main1}, $\text{Out}_{s.g.}(K)$ is isomorphic to $H^1( \mathscr{H}_P, k^\times) \rtimes \text{Aut}(P)$. 
\end{proof}

\section{The Complexity of an Incidence Monoid}\label{S:Complexity}

Let $G$ be a connected linear algebraic group with subgroups $B,U,T$ as before. 
Let $X$ be an irreducible algebraic variety on which $G$ acts morphically. 
As we mentioned in the introduction, the complexity of $X$ with respect to the action of $G$ is the transcendence degree of $\dim k(X)^B$ over $k$. 
It follows from basic principles of algebraic group actions that
\begin{align}\label{A:maindefinitionofc}
c_G(X)= \dim X - \max \{ \dim B\cdot x :\ x\in X\}.
\end{align}
In other words, the complexity of $X$ is the maximal number of parameters determining a continuous family of $B$-orbits in $X$.

\subsection{Proof of Theorem~\ref{T:main2}.}

Recall that Theorem~\ref{T:main2} states that, 
$c_{T\times T}(I(P))$ is given by $Z(P,3)-2Z(P,2)+1$, where $Z(P,n)$ is the zeta polynomial of $P$. 

\begin{proof}[Proof of Theorem~\ref{T:main2}]

Let $P$ be a connected poset with $n$ elements. 
In (\ref{A:maindefinitionofc}), if $G$ is a torus $S$, then we have $G=B=S$, hence, 
$c_S(X)= \dim X -  \max \{ S\cdot x :\ x\in X\}$. 
The irreducible variety $X$ in our case is the incidence monoid $K:= I(P)$, 
and $S$ is the torus $T\times T$ which acts on $I(P)$ by $(s,t)\cdot x = s x t^{-1}$
for every $x\in I(P)$, $s,t\in T$.

We claim that there exists a $T\times T$-orbit of dimension $2n-1$ in $I(P)$. 
Indeed, let $\mathscr{T}$ denote the spanning tree of the Hasse diagram of $P$. 
Since $P$ is connected, $\mathscr{T}$ has $n-1$ edges. 
Let $Q$ denote the poset whose Hasse diagram is $\mathscr{T}$. 
Then $Q$ is a connected subposet of $P$ with at least $n-1$ nontrivial relations. 
In particular, $I(Q)$ is a $T\times T$-stable submonoid of $I(P)$. 
We claim that the kernel of the action of $T\times T$ on $I(Q)$ is the diagonal subgroup,
$\{(c \cdot id_n , c\cdot id_n ):\ c\in k^\times \}$. 
To this end, without loss of generality, we assume that $I(Q)$ is a $k$-subalgebra of $\overline{\mathbf{B}_n}$ 
and we assume that $T= \mathbf{T}_n$. 
Let $\zeta$ denote the zeta function of $I(Q)$, and let 
$(t,s)=(\text{diag}(t_1,\dots,t_n),\text{diag}(s_1,\dots, s_n))$ be an element from the kernel of the action 
of $\mathbf{T}_n\times \mathbf{T}_n$. 
Then for every $x\in K$, we have $txs^{-1} =  x$. 
In particular, if we choose $x$ as the zeta function of the incidence algebra, then we see that 
$t_i = s_j$ for every $x_i \leq x_j$ in $Q$. 
It follows that $t_i = s_i$ for $i\in \{1,\dots, n\}$. 
But since $\mathscr{H}_Q$ is connected, we see that $t_i = s_j$ for all $i,j\in \{1,\dots, n\}$. 
In other words, $(s,t)$ is in the kernel, hence, $\dim \mathbf{T}_n \zeta \mathbf{T}_n = 2n-1$.
Since $\zeta$ is an element of $I(P)$ also, we see that $\max \{ \dim (T\times T)\cdot x :\ x\in I(P) \} = 2n-1$. 
In particular, we proved that 
$c_{T\times T}(I(P))= \dim I(P) - (2n-1) = Z(P,3) - 2Z(P,2) +1$. 
\end{proof}

We have a remark about a notion that is closely related to the complexity. 
\begin{Remark}
Let $X$ be a $G$-variety, where $G$ is a connected linear algebraic group. 
The {\em rank of $G:X$}, denoted by $r_G(X)$, is defined as the rank of the finitely generated abelian group, 
$\varLambda (X):= \{ \chi \in \mathcal{X}(B):\ \text{there exists an $f\in k(X)$ s.t. $b\cdot f = \chi(b) f$ for all $b\in B$}\}$.
The following facts about complexity and rank are well-known (see~\cite[Chapter 2, Section 5]{Timashev}):
\begin{enumerate}
\item  If $X$ is a quasi-affine variety, then $\text{tr.deg}_k k(X)^U = c_G(X) + r_G(X)$. 
\item  If $G$ is a torus $S$, then $c_S(X) = \dim k(X)^S$ and $r_S(X) = \dim S- \dim S_0$, where $S_0$ is the kernel of the action $S:X$.
\end{enumerate}
By combining 1. and 2. we see that $c_S(X) =  \text{tr.deg}_k k(X) - r_S(X) = \dim X - r_S(X)$. 
Thus, it follows from Theorem~\ref{T:main2} that, for $X= I(P)$, and $S=T\times T$, where $T$ is a maximal torus of $I(P)$,
we have $r_{T\times T}(I(P)) = 2|P|-1 = 2Z(P,2)-1$. 
\end{Remark}

\subsection{Which incidence monoids are toric varieties?}\label{SS:Toric}

In this section we will determine all posets $P$ such that $c_{T\times T}(I(P))=0$, where $T$ is a maximal torus of the unit group of $G(P)$. 
As we mentioned in the introduction, these incidence monoids are precisely the toric incidence monoids. 
Recall also that the zeta polynomial of $P$ is the unique polynomial function $Z(P,n)$, 
such that if $n\geq 2$, then $Z(P,n)$ is the number of multichains $t_1 \leq t_2 \leq \cdots \leq t_{n-1}$ in $P$. 
Let $x$ and $y$ be two elements from $P$ such that $x\leq y$. 
If $x=y$, then $x\leq y$ will be called a {\em trivial relation of $P$};
otherwise, $x\leq y$ will be called a {\em nontrivial relation of $P$}. 
Clearly, any covering relation of $P$ is a nontrivial relation of $P$.
Also, it follows from definitions that 
\begin{align}\label{A:Zeta3}
\dim_k I(P) = \text{$\#$ trivial relations + $\#$ nontrivial relations} = Z(P,3). 
\end{align}

\begin{Lemma}\label{L:toricthen}
Let $P$ be a poset with $n$ elements. 
If the incidence monoid $I(P)$ is a toric $\mathbf{T}_n\times \mathbf{T}_n$-variety, then the number of nontrivial relations of $P$ is bounded by $n-1$. 
\end{Lemma}

\begin{proof}
The proof is a consequence of a simple dimension argument:
By (\ref{A:Zeta3}) we know that if $P$ has $m$ nontrivial relations, then $\dim I(P) = n+m$. 
Also, the kernel of the action $\mathbf{T}_n\times \mathbf{T}_n:  I(P)$ contains the diagonal subgroup,
\[
\ker (\mathbf{T}_n\times \mathbf{T}_n:  I(P)) = \{ (c \cdot id_n ,c^{-1}\cdot id_n) :\ c \in k^\times \}. 
\]
Therefore, the dimension of the maximal dimensional orbit is bounded by $2n-1$.
Since $\mathbf{T}_n\times \mathbf{T}_n$ has an open orbit in $I(P)$, we see that $2n-1 \geq n+m$, or $n-1 \geq m$.
This finishes the proof. 
\end{proof}

\begin{Lemma}\label{L:circuits}
Let $Q$ be a connected poset and let $r$ be a positive integer.  
If the Hasse diagram of $Q$ has $r$ circuits, then $\dim I(Q) \geq 2|Q|-1+r$. 
\end{Lemma}

\begin{proof}
Since $\mathscr{H}_Q$ is a connected graph, it has a spanning tree with $|Q|-1$ elements. 
Each edge of this spanning tree corresponds to a covering relation in $Q$. 
To create a circuit, we have to add one edge between two vertices in the spanning tree. 
Therefore, if $\mathscr{H}_Q$ has $r$ circuits, then $Q$ has at least $|Q|-1+r$ nontrivial relations.
The rest of the proof follows from (\ref{A:Zeta3}).
\end{proof}

\begin{Example}
Let $Q$ be the poset whose Hasse diagram is in Figure~\ref{F:twocircuits}.
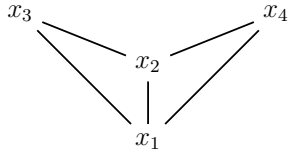
\begin{figure}[htp]
\begin{center}
\scalebox{.85}{
\begin{tikzpicture}[scale=1]

\node (0) at (0,0) {$x_1$};
\node (1) at (0,1.2) {$x_2$};
\node (2) at (-2,2) {$x_3$};
\node (3) at (2,2) {$x_4$};

\draw[thick,-] (0) to (1);
\draw[thick,-] (1) to (2);
\draw[thick,-] (1) to (3);
\draw[thick,-] (0) to (3);
\draw[thick,-] (0) to (2);
\end{tikzpicture}
}
\end{center}
\caption{A poset with 3 circuits.}
\label{F:twocircuits}
\end{figure}
Then $\mathscr{H}_Q$ has 3 circuits and 4 vertices. It is easy to check that $Q$ has 7 nontrivial relations. 
Hence, we have $\dim I(Q) = 11 \geq 8 -1 + 3 = 10$.
Another example is given by the poset in Figure~\ref{F:firstexample}.
Then $\mathscr{H}_P$ has 3 circuits, 5 vertices, and 9 nontrivial relations. 
In this case, we have $\dim I(P) = 14 \geq 10 - 1 + 3 = 12$.

\end{Example}

We are ready to prove the first part of Theorem~\ref{T:main3}, which states that 
if $P$ is a connected poset, then the complexity $c_{T\times T}(I(P))$ is zero if and only if 
all maximal chains of $P$ are of length 1, and the underlying graph of Hasse diagram of $P$ has no circuits.

\begin{proof}[Proof of Part 1 of Theorem~\ref{T:main3}]

Let $P$ be a connected poset with $n$ elements, and let $K$ denote its incidence monoid. 
By fixing a linear extension, $P=\{x_1,\dots, x_n\}$, 
we view $K$ as a $T\times T$-variety, where $T= \mathbf{T}_n$.

We proceed with the assumption that $c_{T\times T}(K) =0$.
Then $K$ is a toric variety, hence, $\dim K \leq 2n-1$.
If the Hasse diagram $\mathscr{H}_P$ has a circuit, then by Lemma~\ref{L:circuits}, $2n < \dim K \leq 2n$,
which is absurd. 
Therefore, $\mathscr{H}_P$ cannot have a circuit.
Let us now assume that $P$ cannot have a saturated chain of length 2. 
Since $\mathscr{H}_P$ has no circuit, this chain is a branch of the spanning tree of $P$. 
We already know that the spanning tree of $\mathscr{H}_P$ gives $n-1$ nontrivial relations. 
If $x \lneq y \lneq z$ is the chain, then we get an additional nontrivial relation by $x \lneq z$. 
Then, in $P$, we have at least $n$ nontrivial relations, which implies that $\dim K \geq 2n$. 
But this contradicts the fact that $K$ is a toric variety of dimension at most $2n-1$. 
Therefore, every maximal chain in $P$ has length 1.

Conversely, let $P$ be a connected poset with $n$ elements such that 
$\mathscr{H}_P$ has no circuits, and every maximal chain in $P$ has length 1. 
Then it is easy to check that $\mathscr{H}_P$ is a bipartite graph with no circuits, hence, $\mathscr{H}_P$ is a tree. 
In particular, it has exactly $n-1$ edges. 
It follows that $\dim K = 2n-1$. 
The rest of the proof follows from (\ref{A:Zeta3}) and Theorem~\ref{T:main2} which states in our case that 
$c_{T\times T}(K) = 0$ if and only if $Z(P,3) = 2n-1$.

\end{proof}

The following corollary is a consequence of the proof of Theorem~\ref{T:main3}.

\begin{Corollary}\label{C:connectedtoricones}
Let $P$ be a connected poset such that $|P|=n$.  
Then the following statements are equivalent:
\begin{enumerate}
\item The incidence monoid $I(P)$ is a toric $\mathbf{T}_n\times \mathbf{T}_n$-variety.
\item The dimension of the incidence monoid $I(P)$ is $2n-1$.
\item The number of nontrivial relations in $P$ is $n-1$.
\item All maximal chains of $P$ have length 1, and the underlying graph of the Hasse diagram of $P$ has no circuits. 
\end{enumerate}
\end{Corollary}

By using posets with multiple connected components, we can get lower dimensional toric varieties. 
The simplest example is when $P$ is an anti-chain. In this case, $I(P)$ is nothing but the monoid of diagonal $n\times n$ matrices. 
Let us prove a converse of Corollary~\ref{C:connectedtoricones} for posets whose connected components have at least one nontrivial relation.

\begin{Proposition}\label{P:rcomponents}
Let $P$ be a poset such that $|P|=n$. Let $P_1,\dots, P_r$ denote the connected components of $P$,
and let $n_1,\dots, n_r$ denote the numbers of nontrivial relations in $P_1,\dots, P_r$, respectively. 
If for every $i\in \{1,\dots, r\}$ we have $n_i = |P_i|-1$, then $I(P)$ is a toric $\mathbf{T}_n\times \mathbf{T}_n$-variety. 
\end{Proposition}

\begin{proof}
The incidence monoid of a disjoint union $P=P_1\sqcup \cdots \sqcup P_r$ is the direct product of incidence algebras, $I(P)=I(P_1)\times \cdots \times I(P_r)$. 
Since $n_i = |P_i|-1$ for every $i\in \{1,\dots, r\}$, by Corollary~\ref{C:connectedtoricones}, every factor $I(P_i)$ is a toric subvariety. The proof follows from the fact that the product of finitely many toric varieties is a toric variety.  
\end{proof}

\begin{Example}
The incidence monoids $I(P)$, where $|P|=3$, such that $c_{T\times T}(I(P)) \leq 1$ are depicted in Figure~\ref{F:n=3}. 
The first four of them are toric incidence monoids.
\begin{figure}[htp]
\begin{center}
\scalebox{.85}{
\begin{tikzpicture}[scale=1]

\begin{scope}[xshift= -8cm]
\node at (0,-1) {$\dim I(P) =3$};
\node (0) at (-1,0) {$x_1$};
\node (1) at (0,0) {$x_2$};
\node (2) at (1,0) {$x_3$};
\end{scope}

\begin{scope}[xshift= -4cm]
\node at (0,-1) {$\dim I(P) =4$};
\node (0) at (0,0) {$x_1$};
\node (1) at (0,1) {$x_2$};
\node (2) at (1,0) {$x_3$};
\draw[thick,-] (0) to (1);
\end{scope}

\begin{scope}[xshift= 0cm]
\node at (0,-1) {$\dim I(P) =5$};
\node (0) at (0,1) {$x_3$};
\node (1) at (-1,0) {$x_1$};
\node (2) at (1,0) {$x_2$};
\draw[thick,-] (0) to (1);
\draw[thick,-] (0) to (2);
\end{scope}

\begin{scope}[xshift= 4cm]
\node at (0,-1) {$\dim I(P) =5$};
\node (0) at (0,0) {$x_1$};
\node (1) at (-1,1) {$x_2$};
\node (2) at (1,1) {$x_3$};
\draw[thick,-] (0) to (1);
\draw[thick,-] (0) to (2);
\end{scope}

\begin{scope}[xshift=8cm]
\node at (0,-1) {$\dim I(P) =6$};
\node (0) at (0,0) {$x_1$};
\node (1) at (0,1) {$x_2$};
\node (2) at (0,2) {$x_3$};
\draw[thick,-] (0) to (1);
\draw[thick,-] (1) to (2);
\end{scope}
\end{tikzpicture}
}
\end{center}
\caption{The incidence monoids of complexity $\leq 1$ with $|P|=3$.}
\label{F:n=3}
\end{figure}
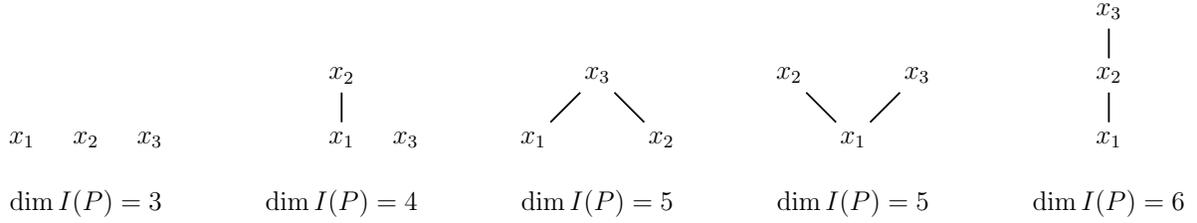

\end{Example}

\begin{Corollary}
Let $n \geq 3$. 
Then for every $r\in \{ n,n+1,\dots,  2n\}$, there is a poset $P=P(r)$ such that $I(P)$ is toric variety of dimension $r$. 
\end{Corollary}
\begin{proof}
The proof follows from Proposition~\ref{P:rcomponents} by induction. 
\end{proof}

We are now ready to prove Part 2. of Theorem~\ref{T:main3}, which states that if $P$ is a connected poset, then 
$c_{T\times T}(I(P))=1$ if and only if $P$ is a graded poset, and one of the following two statements holds:
\begin{enumerate}
\item all maximal chains of $P$ are of length 1, and the underlying graph of the Hasse diagram of $P$ has exactly one circuit; 
\item $P$ has a unique maximal chain of length 2, and the underlying graph of the Hasse diagram of $P$ has no circuits. 
\end{enumerate}

\begin{proof}[Proof of Part 2. of Theorem~\ref{T:main3}]

Let $P$ be a connected poset with $n$ elements, and let $K$ denote its incidence monoid. 
By Theorem~\ref{T:main1}, we know that $c_{T\times T}(K) =1$ if and only if $\dim K = 2n$.

($\Rightarrow$)
Notice that if either 1) or 2) holds, then $P$ is a graded poset. 
Therefore, it suffices to show that one of them holds. 
We proceed with the assumption that 2) does not hold. 
In this case, if $\mathscr{H}_P$ has at least two circuits, then by Lemma~\ref{L:circuits} we obtain
the absurd statement that $2n = 2|P| < \dim K = 2n$. 
Therefore, if 2) does not hold, then 1) holds. 
This finishes the proof of the implication ($\Rightarrow$).

($\Leftarrow$)
Next, we will prove that if $P$ is graded, and either 1) or 2) holds, then $K$ is a toric variety. 
By using the spanning tree argument once more, we see that if $P$ has a chain of length 2, then it cannot have a circuit, and conversely, if $P$ has a circuit, then it cannot have a chain of length 2. 
In other words, the cases 1) and 2) are exclusive conditions. 
Nevertheless, in both of these cases, the number of nontrivial relations in $P$ is exactly $n$, hence, $\dim K = 2n$.
This finishes the proof of our theorem.
\end{proof}

We proceed to compute the automorphism groups of incidence monoids of complexity $\leq 1$.  
For this purpose, the following terminology will be convenient. 
\begin{Definition}\label{D:typesoftoric}
A connected poset $P$ is said to be  
\begin{itemize}
\item a {\em toric poset} if every maximal chain in $P$ has length 1 and 
$\mathscr{H}_P$ has no circuits;
\item a {\em complexity 1 poset with a circuit} if every maximal chain in $P$ has length 1 and 
$\mathscr{H}_P$ contains a single circuit;
\item a {\em complexity 1 poset with a 2-chain} if $P$ is a graded poset such that $P$ has a single 2-chain. 
\end{itemize}
\end{Definition}

\begin{Proposition}
Let $P$ be a connected poset with $n$ elements. 
We assume that $I(P)$ is a toric variety.
Then, the following statements hold: 
\begin{enumerate}
\item If $P$ is a complexity 1 poset with a circuit, then we have $\text{Out}_{s.g.} (I(P)) = k^\times \rtimes \text{Aut}(P)$.
\item If $P$ is either a toric poset, or a complexity 1 poset with a 2-chain, then we have $\text{Out}_{s.g.} (I(P)) = \text{Aut}(P)$.
\end{enumerate}
\end{Proposition}

\begin{proof}
The proof is a direct consequence of Corollary~\ref{C:maincorollary} and our Definition~\ref{D:typesoftoric}.
If $\mathscr{H}_P$ has a circuit, then we have $k^\times$ as a factor since, in this case, we have $H^1(\mathscr{H}_P,k^\times )\cong k^\times$.
Otherwise, we do not have any additional torus factor, and the outer automorphism group of $I(P)$ is given by $\text{Aut}(P)$.
\end{proof}

\begin{Example}
For $n\geq 2$, let $P_n$ denote, as before, the star poset on $n$ vertices.
Then $P_n$ is a toric poset.
The full automorphism group of $P_n$ is given by 
\[
\text{Aut}_{s.g.}(P_n) \cong (G(P_n)/Z) \rtimes S_{n-1},
\] 
where $S_{n-1}$ is the $(n-1)$-th symmetric group, and $Z\cong k^*$ is the center of $G(P_n)$.
\end{Example}

\section{The Incidence Monoid of a Star Poset}\label{S:Star}

The purpose of this section is to give a detailed description of the adherence order of the toric incidence monoid of the star poset on $n$ vertices. 
Since $|P_n| =n$ and the number of nontrivial relations in $P_n$ is $n-1$, 
as an algebraic monoid (or as a $k$-algebra), $I(P_n)$ is $(2n-1)$-dimensional.  
In particular, $I(P_n)$ is a toric variety with respect to $T\times T$, where $T$ is a maximal torus in $I(P_n)$.

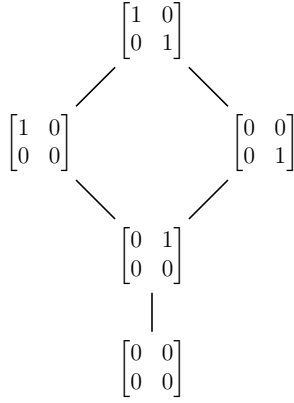
\begin{figure}[htp]
\begin{center}
\scalebox{.75}{
\begin{tikzpicture}[scale=1]
\node (1) at (0,-3) {$\begin{bmatrix} 0 & 0 \\ 0 & 0 \end{bmatrix}$};
\node (2) at (0,-1) {$\begin{bmatrix} 0 & 1 \\ 0 & 0 \end{bmatrix}$};
\node (3) at (-2,1) {$\begin{bmatrix} 1 & 0 \\ 0 & 0 \end{bmatrix}$};
\node (4) at (2,1) {$\begin{bmatrix} 0 & 0 \\ 0 & 1 \end{bmatrix}$};
\node (5) at (0,3) {$\begin{bmatrix} 1 & 0 \\ 0 & 1 \end{bmatrix}$};
\draw[thick,-] (1) to (2); 
\draw[thick,-] (2) to (3);
\draw[thick,-] (2) to (4); 
\draw[thick,-] (3) to (5); 
\draw[thick,-] (4) to (5);  
\end{tikzpicture}
}
\end{center}
\caption{The adherence order on $I(P_2)$.}
\label{F:BCR2}
\end{figure}

For $n=2$, the incidence monoid $I(P_2)$ is nothing but the Borel monoid of all $2\times 2$ upper triangular matrices. 
In this case, the adherence order is given by the Bruhat-Chevalley-Renner order on the interval $[0,id_2]$
in the rook monoid $R_2$. 
The Hasse diagram of the corresponding poset is shown in Figure~\ref{F:BCR2}.

For $n\geq 3$, the elements of $I(P_n)$ are of the form, 
\begin{align}\label{A:xasin}
x=\begin{bmatrix}
x_{11} & x_{12} & x_{13} & \cdots & x_{1n} \\
0 & x_{22} & 0 & \cdots & 0 \\
0 & 0 & x_{33} & \cdots & 0 \\
\vdots & \vdots & \vdots & \ddots & \vdots \\
0 & 0 & 0 & \cdots & x_{nn}
\end{bmatrix}.
\end{align} 
For a given $x$ as in (\ref{A:xasin}), we want to understand the set of all products of the form $axb$, 
where $a$ and $b$ are invertible elements from $I(P_n)$.
Let $a$ and $b$ be given by 
\begin{align}\label{A:AandB}
a=\begin{bmatrix}
a_{11} & a_{12} & a_{13} & \cdots & a_{1n} \\
0 & a_{22} & 0 & \cdots & 0 \\
0 & 0 & a_{33} & \cdots & 0 \\
\vdots & \vdots & \vdots & \ddots & \vdots \\
0 & 0 & 0 & \cdots & a_{nn}
\end{bmatrix}\ \text{ and }\
b=\begin{bmatrix}
b_{11} & b_{12} & b_{13} & \cdots & b_{1n} \\
0 & b_{22} & 0 & \cdots & 0 \\
0 & 0 & b_{33} & \cdots & 0 \\
\vdots & \vdots & \vdots & \ddots & \vdots \\
0 & 0 & 0 & \cdots & b_{nn}
\end{bmatrix}.
\end{align}
Let $y_{ij}= y_{ij}(a,x,b)$ denote the $(i,j)$-th entry of $axb$. 
Then we have 
\begin{align}\label{A:ijthentry}
y_{ij} =
\begin{cases}
a_{ii}x_{ii}b_{ii} & \text{ if $i=j$};\\
0 & \text{ if $i\neq j$ and $i\neq 1$};\\
a_{11} x_{11} b_{1j} + a_{11}x_{1j}b_{jj} + a_{1j}x_{jj}b_{jj}& \text{ if $i\neq j$ and $i=1$}.
\end{cases}
\end{align}

Now we are going to investigate the structure of the inclusion poset of $G(P_n)\times G(P_n)$-orbit closures.

\begin{Lemma}\label{L:nox11}
Let $O$ be a $G(P_n)\times G(P_n)$-orbit in $I(P_n)$.
Let $x = (x_{i,j})_{i,j=1}^n$ be an element from $O$. 
If $x_{11}\neq0$, then $O$ contains a diagonal rook matrix. 
\end{Lemma}

\begin{proof}
Let $a$ and $b$ be two matrices from $G(P_n)$ whose entries are defined as in (\ref{A:AandB}).
Let $(y_{ij})_{i,j=1}^n$ denote the matrix $axb$. 
As we observed in (\ref{A:ijthentry}), the $(1,j)$-th entry of $(y_{ij})_{i,j=1}^n$ is given by 
$y_{1j} = a_{11} x_{11} b_{1j} + a_{11}x_{1j}b_{jj} + a_{1j}x_{jj}b_{jj}$.
Since $a$ and $b$ are from $G(P_n)$, we know that $a_{jj},b_{jj} \in k^*$ and 
$a_{1j},b_{1j}\in k$ for every $j\in \{1,\dots, n\}$.
Then, for each $j\in \{2,\dots, n\}$, we choose $b_{1j}$ so that the following equality holds, 
\[
-a_{11} x_{11} b_{1j} = a_{11}x_{1j}b_{jj} + a_{1j}x_{jj}b_{jj}.
\]
In other words, by choosing $b\in G(P_n)$ appropriately, we find a matrix $(y_{ij})_{i,j=1}^n = a x b$ in $G(P_n)x G(P_n)$
such that $y_{1j} =0$ for every $j\in \{2,\dots, n\}$. 
It follows that the orbit $O$ is represented by a diagonal matrix. 
Then, by multiplying with a suitable (invertible) diagonal matrix, we make the entries of this representative either 0 or 1;
hence, $O$ has a diagonal representative which is a rook matrix. 
This finishes the proof of our lemma.
\end{proof}

Next, we will prove the complimentary result.

\begin{Lemma}\label{L:nox1j}
Let $O$ be a $G(P_n)\times G(P_n)$-orbit in $I(P_n)$.
Let $x = (x_{i,j})_{i,j=1}^n$ be an element from $O$. 
If $x_{11}=0$, then $O$ contains a 0/1 matrix $(y_{ij})_{i,j=1}^n$ such that 
\begin{enumerate}
\item if $x_{jj} \neq 0$ for some $j\in \{2,\dots, n\}$, then $y_{1j}=0$ and $y_{jj}=1$;
\item if $x_{jj}=0$, and $x_{1j}\neq 0$ for some $j\in \{2,\dots, n\}$, then $y_{1j}=1$. 
\end{enumerate}
\end{Lemma}

\begin{proof}
As in the proof of previous lemma, 
we let $(y_{ij})_{i,j=1}^n$ denote the matrix $axb$, where  $a$ and $b$ are as in (\ref{A:AandB}).
We fix the entry $a_{11}$ of $a$. 

Now we assume that $x_{jj}\neq 0$ for some $j>1$. 
By (\ref{A:ijthentry}), the $(1,j)$-th entry of $(y_{ij})_{i,j=1}^n$ is given by 
$y_{1j} = a_{11}x_{1j}b_{jj} + a_{1j}x_{jj}b_{jj}$. 
We fix the entry $b_{jj} \in k^\times$ and we chose $a_{1j}$ so that $a_{11}x_{1j}b_{jj} + a_{1j}x_{jj}b_{jj} = 0$ holds. 
Then, we have $y_{1j}=0$. 
Likewise, by (\ref{A:ijthentry}), we know that the $(j,j)$-th entry of $(y_{ij})_{i,j=1}^n$ is given by $a_{jj}x_{jj} b_{jj}$. 
Then we let $a_{jj} = \frac{1}{x_{jj}b_{jj}}$, hence, we get $y_{jj}=1$. 

Next, we assume that $x_{jj}= 0$, and $x_{1j}\neq 0$ for some $j>1$. 
By (\ref{A:ijthentry}), the $(1,j)$-th entry of $(y_{ij})_{i,j=1}^n$ is given by 
$y_{1j} = a_{11}x_{1j}b_{jj}$. 
Then we chose $b_{jj}$ so that $a_{11}x_{1j}b_{jj}  = 1$ holds; hence, we find that $y_{1j}=1$. 
This finishes the proof of our lemma.
\end{proof}

We summarize Lemmas~\ref{L:nox11} and~\ref{L:nox1j}.

\begin{Proposition}\label{P:elements1}
Let $O$ be a $G(P_n)\times G(P_n)$-orbit in $I(P_n)$. 
Then $O$ is uniquely represented by a 0/1 matrix $(y_{ij})_{i,j=1}^n$ contained in $I(P_n)$ such that 
\begin{enumerate}
\item there exists at most one 1 in each column of $(y_{ij})_{i,j=1}^n$;
\item if $y_{11} = 1$, then $y_{1j}=0$ for every $j\in \{2,\dots, n\}$.
\end{enumerate}
\end{Proposition}

\begin{proof}
In light of Lemmas~\ref{L:nox11} and~\ref{L:nox1j}, we have to show the uniqueness of the representatives only. 
First, if $y_{11}\neq 0$, then we have a diagonal rook matrix as the representative.
Since diagonal rook matrices are the representatives for $\mathbf{T}_n\times \mathbf{T}_n$-orbits in $\overline{\mathbf{T}_n}$,
$y$ is unique.  
Next, if $y_{11}=0$ holds, then once again by Lemma~\ref{L:nox1j}, 
we have an orbit representative which is a 0/1 matrix. 
If two such matrices are in the same $G(P_n)\times G(P_n)$-orbit, then first observe that their diagonal entries are the same.
After this, it is easy to see that the first rows of two such matrices are equal.
Hence, we have a unique orbit representative when $y_{11}=0$, also. 
This finishes the proof of our assertion. 
\end{proof}

\begin{Notation}
Hereafter we will denote the set of $G(P_n)\times G(P_n)$-orbit representatives of $I(P_n)$ that we described in
Proposition~\ref{P:elements1} by $R_{P_n}$.
\end{Notation}

\begin{Corollary}
Let $n\geq 2$. If $x$ is an element from $R_{P_n}$, then $x$ is a partial transformation of the set $\{1,2,\dots, n\}$. 
\end{Corollary}
\begin{proof}
Since every column of $x$ has at most one nonzero entry, the proof follows. 
\end{proof}

\begin{Example}\label{E:elementsofRP3}
The elements of $R_{P_3}$ are the following: 
\begin{align*}
\begin{bmatrix}
0 & 0 & 0 \\
0 & 0 & 0 \\
0& 0 & 0
\end{bmatrix},
\begin{bmatrix}
1 & 0 & 0 \\
0 & 0 & 0 \\
0& 0 & 0
\end{bmatrix},
\begin{bmatrix}
0 & 0 & 0 \\
0 & 1 & 0 \\
0& 0 & 0
\end{bmatrix},
\begin{bmatrix}
0 & 0 & 0 \\
0 & 0 & 0 \\
0& 0 & 1
\end{bmatrix},
\begin{bmatrix}
1 & 0 & 0 \\
0 & 1 & 0 \\
0& 0 & 0
\end{bmatrix},
\begin{bmatrix}
1 & 0 & 0 \\
0 & 0 & 0 \\
0& 0 & 1
\end{bmatrix},
\begin{bmatrix}
0 & 0 & 0 \\
0 & 1 & 0 \\
0& 0 & 1
\end{bmatrix},\\
\begin{bmatrix}
1 & 0 & 0 \\
0 & 1 & 0 \\
0& 0 & 1
\end{bmatrix},
\begin{bmatrix}
0 & 1 & 0 \\
0 & 0 & 0 \\
0& 0 & 0
\end{bmatrix},
\begin{bmatrix}
0 & 0 & 1 \\
0 & 0 & 0 \\
0& 0 & 0
\end{bmatrix},
\begin{bmatrix}
0 & 1 & 1 \\
0 & 0 & 0 \\
0& 0 & 0
\end{bmatrix},
\begin{bmatrix}
0 & 1 & 0 \\
0 & 0 & 0 \\
0& 0 & 1
\end{bmatrix},
\begin{bmatrix}
0 & 0 & 1 \\
0 & 1 & 0 \\
0& 0 & 0
\end{bmatrix}.
\end{align*}

\end{Example}

By using the concrete description of the elements of $R_{P_n}$, we are able to determine the size of $R_{P_n}$. 

\begin{Lemma}
Let $n \geq 2$. 
The number of $G(P_n)\times G(P_n)$-orbits in $I(P_n)$ is $2^{n-1}+3^{n-1}$. 
\end{Lemma}
\begin{proof}
First, we count the elements of $R_{P_n}$ that have 1 at their $(1,1)$-th entry.
We know from Lemma~\ref{L:nox11} that such representatives are diagonal rook matrices. 
Since we already know that the upper left corner gets a 1, the remaining entries of such a representative can be chosen in $2^{n-1}$ different ways. 
Next, we count the number of $y\in R_{P_n}$ such that $y= (y_{ij})_{i,j=1}^n$ and $y_{11}=0$. 
In this case, we choose $k$ diagonal entries, and use the remaining $n-1-k$ columns and the first row 
to place 1's. The total number of choices is $2^{n-1-k}$.  
In summary, we have the following formula for the cardinality of $R_{P_n}$:
\begin{align}\label{A:formula}
|R_{P_n}| &= 2^{n-1} + \sum_{k=0}^{n-1} {n-1 \choose k} 2^{n-1-k}.
\end{align}
It is easy to show that the summation on the right hand side of (\ref{A:formula}) is equal to $3^{n-1}$. 
This finishes the proof of our lemma.
\end{proof}

\begin{Lemma}\label{L:dimension}
Let $x:=(x_{ij})_{i,j=1}^n$ be an element from $R_{P_n}$. 
Then the dimension of the corresponding orbit is given by 
\[
\dim G(P_n) x G(P_n) = 
\begin{cases}
n-1 + \sum_{i=1}^n x_{ii} & \text{ if $x_{11}=1$};\\
\sum_{i=2}^n 2x_{ii} + \sum_{i=2}^n x_{1i}   & \text{ if $x_{11}=0$}.
\end{cases}
\]
\end{Lemma}

\begin{proof}
Let us first assume that $x_{11}=1$. 
Then we see from~(\ref{A:ijthentry}) that for every $r\in \{2,\dots, n\}$, $y_{1r} \neq 0$ and $y_{rr}\neq 0$ for the appropriate choices of $a$ and $b$.
It follows that, if $x_{11}=1$, then we have exactly $n-1 + \sum_{i=1}^n x_{ii}$ free parameters for the orbit $G(P_n)xG(P_n)$.

Next, we assume that $x_{11}=0$. 
In this case, we see from~(\ref{A:ijthentry}) that for every $r\in \{2,\dots, n\}$ such that $x_{rr}\neq 0$,
we have exactly two free parameters: $y_{1r}$ and $y_{rr}$. 
At the same time, for every $r\in \{2,\dots, n\}$ such that $x_{1r}\neq 0$ we have exactly one free parameter.  
Since these two cases are disjoint, the total number of free parameters for $G(P_n)xG(P_n)$ is 
given by $\sum_{i=2}^n 2x_{ii} + \sum_{i=2}^n x_{1i}$.  
This finishes the proof of our assertion. 
\end{proof}

We are now ready to determine the covering relations of $R_{P_n}$.
To ease our notation, if $n$ is a positive integer, then we will denote the set $\{1,\dots, n\}$ by $[n]$.

\begin{Theorem}\label{T:adherence}
Let  
\begin{align*}
y=\begin{bmatrix}
y_{11} & y_{12} & y_{13} & \cdots & y_{1n} \\
0 & y_{22} & 0 & \cdots & 0 \\
0 & 0 & y_{33} & \cdots & 0 \\
\vdots & \vdots & \vdots & \ddots & \vdots \\
0 & 0 & 0 & \cdots & y_{nn}
\end{bmatrix}\ \text{ and }\
x=\begin{bmatrix}
x_{11} & x_{12} & x_{13} & \cdots & x_{1n} \\
0 & x_{22} & 0 & \cdots & 0 \\
0 & 0 & x_{33} & \cdots & 0 \\
\vdots & \vdots & \vdots & \ddots & \vdots \\
0 & 0 & 0 & \cdots & x_{nn}
\end{bmatrix}
\end{align*}
be two distinct elements from $R_{P_n}$. 
Then $y$ is covered by $x$ in the adherence order if and only if one of the following holds:
\begin{enumerate}

\item There exists $r\in [n]$ such that for every $(i,j)\in [n]\times [n]\setminus \{(r,r)\}$ 
we have 
\begin{enumerate}
\item $x_{ij}=y_{ij}$, 
\item $x_{11}=y_{11}=1$, 
\item $x_{rr}=1$ and $y_{rr}=0$.
\end{enumerate}

\item There exists $r\in \{2,\dots, n\}$ such that for every $(i,j)\in [n]\times [n]\setminus \{(1,r),(1,1)\}$ 
we have 
\begin{enumerate}
\item $x_{ij}=y_{ij}$, 
\item $x_{11}=1$ and $y_{11}=0$, 
\item $x_{1r}=0$, $y_{1r}=1$.
\end{enumerate}

\item $y$ and $x$ are given by 
\begin{align*}
y=\begin{bmatrix}
0 & 1 & 1 & \cdots & 1 \\
0 & 0 & 0 & \cdots & 0 \\
0 & 0 & 0 & \cdots & 0 \\
\vdots & \vdots & \vdots & \ddots & \vdots \\
0 & 0 & 0 & \cdots & 0
\end{bmatrix}\ \text{ and }\
x=\begin{bmatrix}
1 & 0 & 0 & \cdots & 0 \\
0 & 0 & 0 & \cdots & 0 \\
0 & 0 & 0 & \cdots & 0 \\
\vdots & \vdots & \vdots & \ddots & \vdots \\
0 & 0 & 0 & \cdots & 0
\end{bmatrix}.
\end{align*}

\item There exists $r\in \{2,\dots, n\}$ such that for every $(i,j)\in [n]\times [n]\setminus \{(1,r),(r,r)\}$ 
we have 
\begin{enumerate}
\item $x_{ij}=y_{ij}$, 
\item $x_{11}=0$ and $y_{11}=0$, 
\item $x_{rr}=1$ (hence, $x_{1r}=0$), $y_{1r}=1$ (hence, $y_{rr}=0$).
\end{enumerate}

\item There exists $r\in \{2,\dots, n\}$ such that for every $(i,j)\in [n]\times [n]\setminus \{(1,r),(r,r)\}$ 
we have 
\begin{enumerate}
\item $x_{ij}=y_{ij}$, 
\item $x_{11}=0$ and $y_{11}=0$, 
\item $x_{1r}=1$ (hence, $x_{rr}=0$), $y_{1r}=0$.
\end{enumerate}

\end{enumerate}
\end{Theorem}

\begin{proof}

Let $x$ be an element of $R_{P_n}$. 
We will determine all $y\in R_{P_n}$ such that $y$ lies in the closure $\overline{G(P_n)x G(P_n)}$, 
and $\dim G(P_n) y G(P_n) + 1 = \dim G(P_n)x G(P_n)$ holds. 
Now let $a$ and $b$ be two generic elements from $G(P_n)$ such that $a:=(a_{ij})_{i,j=1}^n,b=(b_{ij})_{i,j=1}^n$.

Let us begin with the case where $x$ is as in 3).
In this case, by~(\ref{A:ijthentry}), the orbit $G(P_n)x G(P_n)$ consists of matrices of the form 
\begin{align}\label{A:limit1}
\begin{bmatrix}
a_{11}b_{11} & a_{11} b_{12} &  a_{11} b_{13}  & \cdots &  a_{11} b_{1n}  \\
0 & 0 & 0 & \cdots & 0 \\
0 & 0 & 0 & \cdots & 0 \\
\vdots & \vdots & \vdots & \ddots & \vdots \\
0 & 0 & 0 & \cdots & 0
\end{bmatrix}.
\end{align}
Likewise, if $y$ is as in 3), then by~(\ref{A:ijthentry}), the orbit $G(P_n)y G(P_n)$ consists of matrices of the form 
\begin{align}\label{A:limit2}
\begin{bmatrix}
0 & a_{11} b_{22}  &  a_{11} b_{33}  & \cdots &  a_{11} b_{nn} \\
0 & 0 & 0 & \cdots & 0 \\
0 & 0 & 0 & \cdots & 0 \\
\vdots & \vdots & \vdots & \ddots & \vdots \\
0 & 0 & 0 & \cdots & 0
\end{bmatrix}.
\end{align}
By taking the limits of the matrices in (\ref{A:limit1}) as $b_{11}\to 0$, we get every matrix of the form (\ref{A:limit2}).
This argument shows that $G(P_n)y G(P_n) \subseteq \overline{G(P_n)x G(P_n)}$.
Since $G(P_n)y G(P_n)$ is one codimensional in $\overline{G(P_n)x G(P_n)}$, 
we see not only that $y$ is covered by $x$ but also that it is the unique element in $R_{P_n}$ which is covered by $x$.

Next, let us assume that $x$ has at least two nonzero entries, and that $x_{11}=1$.
In this case, by Lemma~\ref{L:nox11}, we know that $x_{1r}=0$ for every $r\in \{2,\dots, n\}$. 
By Lemma~\ref{L:dimension}, we have 
\[
\dim G(P_n) x G(P_n) = n-1 + \# \text{ nonzero diagonal entries of $x$}.
\]
Now, let $y= (y_{ij})_{i,j=1}^n$ be an element of $R_{P_n}$ such that 
\begin{align*}
G(P_n) y G(P_n) \subseteq \overline{ G(P_n) x G(P_n)}\ \text{ and } \ 
\dim G(P_n) x G(P_n) = \dim G(P_n) y G(P_n)+1.
\end{align*}
If $y_{11}=1$, then $y$ and $x$ differ from each other at exactly one diagonal entry. 
In particular, in this case, we get a covering relation as in 1). 
If $y_{11}=0$, then to get the one codimensional orbit $G(P_n) y G(P_n)$ in $\overline{G(P_n) x G(P_n)}$, 
for every $r\in \{2,\dots, n\}$ such that the $r$-th column of $x$ is 0, we must have $y_{1r}=1$. 
In this case, we get the covering relation that is described in 2).

Next, we proceed with the assumption that $x$ has at least two nonzero entries, and $x_{11}= 0$.
On one hand, if $x_{ll}=1$ for some $l\in \{2,\dots, n\}$, then $x_{1l}=0$ by Lemma~\ref{L:nox1j}. 
By using (\ref{A:ijthentry}), we see that $y_{ll} = a_{ll}  b_{ll}$ and $y_{1l}=  a_{1l} b_{ll}$. 
Let us set $b$ as $id_n$. Also, let us set $a_{1l}:=1$ and $a_{rs}:=0$ for every 
$(r,s)\in [n]\times [n] \setminus \{(1,l)\}$ such that $r\neq s$. 
Now we define $y:=(y_{ij})_{i,j=1}^n$ by $y=\lim_{\substack{a_{ss}\to 1\\ s\neq l}} ( \lim_{ a_{ll}\to 0 } ax)$.
Then, first of all, $y \in \overline{G(P_n)x G(P_n)}$.
Secondly, we have $x_{ij}=y_{ij}$ for $(i,j) \notin \{ (1,l), (l,l)\}$, and $y_{1l} = 1$.
Finally, notice that the orbit $\dim G(P_n)y  G(P_n)$ is one codimensional in $\overline{G(P_n)x G(P_n)}$.
This shows that we obtained a covering relation of type 4). 
On the other hand, if $x_{ll}=0$ and $x_{1l}=1$ for some $l\in \{2,\dots, n\}$, then $x_{11}=0$ and $x_{ll}=0$ by Lemma~\ref{L:nox1j}. 
Once more, by using (\ref{A:ijthentry}), we know that $y_{1l}= a_{ll} b_{ll}$. 
In this case, we set $b:=id_n$, and we let $a$ be of the form 
$\text{diag}(a_{11},\dots, a_{nn})$.
Then, as before, we define $y$ by the limit 
$y:=\lim_{\substack{a_{ss}\to 1\\ s\neq l}} ( \lim_{ a_{ll}\to 0 } ax)$. 
The resulting matrix $y:=(y_{ij})_{i,j=1}^n$ is contained in $\overline{G(P_n)x G(P_n)}$.
Furthermore, we have $x_{ij}= y_{ij}$ for $(i,j) \neq (1,l)$, $x_{1l}=1$, and $y_{1l} = 0$.
As before, we have $\dim G(P_n)x G(P_n) = \dim G(P_n)y  G(P_n)+1$. 
In other words, we obtained a covering relation of type 5). 
Thus, we finished checking all possible covering relations in $R_{P_n}$.
\end{proof}

\begin{Example}
We continue with our Example~\ref{E:elementsofRP3}. 
Recall that the elements of $R_{P_n}$ are partial transformations. 
For $n=3$, let us write the elements of $R_{P_n}$ using word notation in the order that we listed them in Example~\ref{E:elementsofRP3}:
\[
000, 100,020,003,120,103,023,123,010,001,011,013,021.
\]
Notice that the only non-rook element in this list is $011$. 
Now, by using Theorem~\ref{T:adherence},
we depict the Hasse diagram of the adherence order on $R_{P_3}$ in Figure~\ref{F:inclusion}. 

\begin{figure}[htp]
\begin{center}
\scalebox{1}{
\begin{tikzpicture}[scale=1.5]

\node (0) at (0,0) {$000$};
\node (1) at (-1,1) {$010$};
\node (2) at (1,1) {$001$};
\node (3) at (0,2) {$011$};
\node (4) at (-2,2) {$020$};
\node (5) at (2,2) {$003$};

\node (6) at (-2,3) {$021$};
\node (7) at (0,3) {$100$};
\node (8) at (2,3) {$013$};

\node (9) at (-2,4) {$120$};
\node (10) at (0,4) {$023$};
\node (11) at (2,4) {$103$};

\node (12) at (0,5) {$123$};

\draw[thick,-] (0) to (1);
\draw[thick,-] (0) to (2);
\draw[thick,-] (1) to (3);
\draw[thick,-] (1) to (4);
\draw[thick,-] (2) to (3);
\draw[thick,-] (2) to (5);

\draw[thick,-] (3) to (6);
\draw[thick,-] (4) to (6);
\draw[thick,-] (3) to (7);
\draw[thick,-] (3) to (8);
\draw[thick,-] (5) to (8);

\draw[thick,-] (6) to (9);
\draw[thick,-] (6) to (10);

\draw[thick,-] (7) to (9);
\draw[thick,-] (7) to (11);

\draw[thick,-] (8) to (10);
\draw[thick,-] (8) to (11);

\draw[thick,-] (10) to (12);
\draw[thick,-] (11) to (12);
\draw[thick,-] (9) to (12);

\end{tikzpicture}
}

\end{center}

\caption{The adherence order on $R_{P_3}$.}
\label{F:inclusion}
\end{figure}
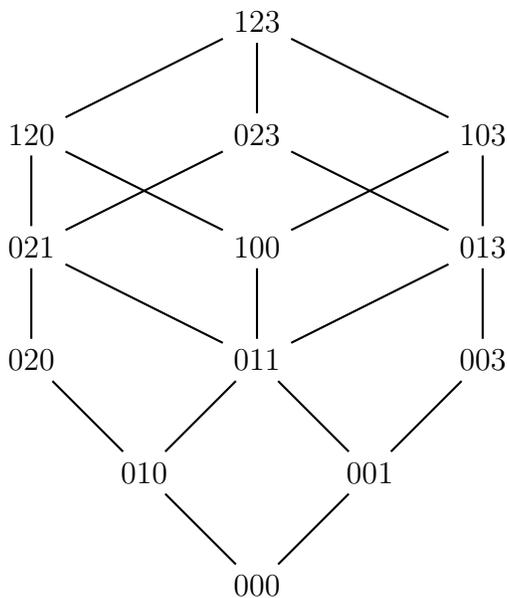
\end{Example}

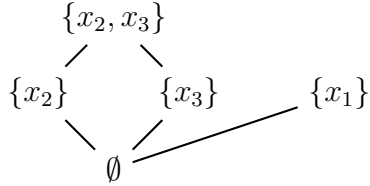
\begin{figure}[htp]
\begin{center}
\scalebox{1}{
\begin{tikzpicture}[scale=1]
\node (0) at (0,0) {$\emptyset$};
\node (1) at (-1,1) {$\{x_2\}$};
\node (2) at (1,1) {$\{x_3\}$};
\node (3) at (0,2) {$\{x_2,x_3\}$};
\node (4) at (3,1) {$\{x_1\}$};
\draw[thick,-] (0) to (1);
\draw[thick,-] (0) to (2);
\draw[thick,-] (1) to (3);
\draw[thick,-] (2) to (3);
\draw[thick,-] (0) to (4);
\end{tikzpicture}
}
\end{center}
\caption{The intersection lattice of antichains for $n=3$.}
\label{F:intersectionlattice}
\end{figure}

Before we present the proof of our Theorem~\ref{T:main4}, let us compute the intersection lattice of the smooth, irreducible, completely regular algebraic submonoids with unipotent radical $U_{P_n}$ of $I(P_n)$ that we introduced in Subsection~\ref{SS:LAM}. 
We continue to use the linear extension of $P_n$ that is defined in Figure~\ref{F:toric1}.
A subset $A$ of $P_n$ is an antichain if either $A \subseteq \{ x_2,\dots, x_n\}$ or $A=\{x_1\}$. 
Then the Hasse diagram of the intersection lattice of the set of antichains is given by 
the Hasse diagram of the Boolean algebra $B_{n-1}$ with an additional vertex which covers the minimum and not related to any other elements. 
For example, for $n=3$, the intersection lattice is as depicted in Figure~\ref{F:intersectionlattice}.
The corresponding completely regular submonoids are given by 
\begin{align*}
I(P_3)_{\{x_2,x_3\}} &= 
\left\{
\begin{bmatrix}
1 & a_1 & a_2 \\
0 & t_2 & 0 \\
0 & 0 & t_3
\end{bmatrix} :\ a_1,a_2,t_2,t_3\in k 
\right\},\\
I(P_3)_{\{x_2\}} &= 
\left\{
\begin{bmatrix}
1 & a_1 & a_2 \\
0 & t_2 & 0 \\
0 & 0 & 1
\end{bmatrix} :\ a_1,a_2, t_2\in k
\right\},\\
I(P_3)_{\{x_3\}} &= 
\left\{
\begin{bmatrix}
1 & a_1 & a_2 \\
0 & 1 & 0 \\
0 & 0 & t_3
\end{bmatrix} :\ a_1,a_2,t_3\in k
\right\},\\
I(P_3)_{\{x_1\}} &= 
\left\{
\begin{bmatrix}
t_1 & a_1 & a_2 \\
0 & 1 & 0 \\
0 & 0 & 1
\end{bmatrix} :\ a_1,a_2, t_1\in k 
\right\},\\
I(P_3)_{\emptyset} &= 
\left\{
\begin{bmatrix}
1 & a_1 & a_2 \\
0 & 1 & 0 \\
0 & 0 & 1
\end{bmatrix} :\ a_1,a_2\in k \right\}.
\end{align*}

We are now ready to finish our paper by presenting the proof of the last main result of our paper. 
\begin{proof}[Proof of Theorem~\ref{T:main4}.]
We already proved parts 1., 2., and 4. of Theorem~\ref{T:main4} earlier in this subsection.
Our third assertion, which states that the adherence order for the action $T\times T:I(P)$ 
is a Boolean algebra of size $2^{n-1}$ needs a proof.  
Once again, we work with the maximal torus $T=\mathbf{T}_n$.
Let $x= (x_{ij})_{i,j=1}^n$ be an element of $I(P_n)$. 
Let $a= \text{diag}(a_1,\dots, a_n)$, $b= \text{diag}(b_1,\dots, b_n)$ be two elements from $T$.
Then we have 
\begin{align}\label{A:tor}
axb=
\begin{bmatrix}
a_{11}x_{11} b_{11} & a_{11}x_{12} b_{22}& a_{11}x_{13} b_{33} & \cdots & a_{11}x_{1n} b_{nn} \\
0 & a_{22}x_{22} b_{22} & 0 & \cdots & 0  \\
0 & 0 & a_{33}x_{33} b_{33} & \dots & 0  \\
\vdots & \vdots & \vdots &  \ddots  &\vdots \\
0 & 0 & \dots & & a_{nn} x_{nn} b_{nn} 
\end{bmatrix}.
\end{align}
Let $A$ denote the set of relevant pairs of indices, $A:= \{ (1,1),\dots, (n,n), (1,2),\dots, (1,n)\}$. 
It follows from (\ref{A:tor}) that the $T\times T$-orbit of $x$ is uniquely determined by the 
set of nonzero entries of $x$.  
In other words, the subset $Z(x):= \{ (i,j) \in A:\ x_{ij}=0\}$ uniquely determines the orbit $TxT$.
Therefore, the total number of $T\times T$-orbits is given by $2^{\dim I(P)}=2^{2n-1}$. 
It is also easy to check from (\ref{A:tor}) that, for $x,y\in I(P)$, 
$TyY\subseteq \overline{TxT}$ if and only if $Z(x) \subseteq Z(y)$. 
Hence, the adherence order of $T\times T$-orbits in $I(P_n)$ is isomorphic to the Boolean 
algebra of subsets of $A$. This finishes the proof of our theorem. 
\end{proof}

\section*{Acknowledgement}

We are grateful to the referee whose suggestions improved the quality of our paper.

\bibliographystyle{plain}
\bibliography{referenc}

\begin{thebibliography}{10}

\bibitem{Backlawski1972}
Kenneth Baclawski.
\newblock Automorphisms and derivations of incidence algebras.
\newblock {\em Proc. Amer. Math. Soc.}, 36:351--356, 1972.

\bibitem{BrusamarelloLewis}
Rosali Brusamarello and David~W. Lewis.
\newblock Automorphisms and involutions on incidence algebras.
\newblock {\em Linear Multilinear Algebra}, 59(11):1247--1267, 2011.

\bibitem{Can2019}
Mahir~Bilen Can.
\newblock The rook monoid is lexicographically shellable.
\newblock {\em European J. Combin.}, 81:265--275, 2019.

\bibitem{CC}
Mahir~Bilen Can and Yonah Cherniavsky.
\newblock Stirling posets.
\newblock {\em Israel J. Math.}, 237(1):185--219, 2020.

\bibitem{CanRenner}
Mahir~Bilen Can and Lex~E. Renner.
\newblock Bruhat-{C}hevalley order on the rook monoid.
\newblock {\em Turkish J. Math.}, 36(4):499--519, 2012.

\bibitem{DrozdKolesnik}
Yuriy Drozd and Petro Kolesnik.
\newblock Automorphisms of incidence algebras.
\newblock {\em Comm. Algebra}, 35(12):3851--3854, 2007.

\bibitem{Grosshans}
Frank~D. Grosshans.
\newblock {\em Algebraic homogeneous spaces and invariant theory}, volume 1673
  of {\em Lecture Notes in Mathematics}.
\newblock Springer-Verlag, Berlin, 1997.

\bibitem{Huang2001}
Wenxue Huang.
\newblock Reductive and semisimple algebraic monoids.
\newblock {\em Forum Math.}, 13(4):495--504, 2001.

\bibitem{Okninski2014}
Jan Okni\'{n}ski.
\newblock On certain semigroups derived from associative algebras.
\newblock In {\em Algebraic monoids, group embeddings, and algebraic
  combinatorics}, volume~71 of {\em Fields Inst. Commun.}, pages 233--245.
  Springer, New York, 2014.

\bibitem{PPR}
Edwin~A. Pennell, Mohan~S. Putcha, and Lex~E. Renner.
\newblock Analogue of the {B}ruhat-{C}hevalley order for reductive monoids.
\newblock {\em J. Algebra}, 196(2):339--368, 1997.

\bibitem{Putcha1983}
Mohan~S. Putcha.
\newblock On the automorphism group of a linear algebraic monoid.
\newblock {\em Proc. Amer. Math. Soc.}, 88(2):224--226, 1983.

\bibitem{Putcha}
Mohan~S. Putcha.
\newblock {\em Linear algebraic monoids}, volume 133 of {\em London
  Mathematical Society Lecture Note Series}.
\newblock Cambridge University Press, Cambridge, 1988.

\bibitem{Renner1985}
Lex~E. Renner.
\newblock Classification of semisimple algebraic monoids.
\newblock {\em Trans. Amer. Math. Soc.}, 292(1):193--223, 1985.

\bibitem{Renner1989}
Lex~E. Renner.
\newblock Completely regular algebraic monoids.
\newblock {\em J. Pure Appl. Algebra}, 59(3):291--298, 1989.

\bibitem{Renner}
Lex~E. Renner.
\newblock {\em Linear algebraic monoids}, volume 134 of {\em Encyclopaedia of
  Mathematical Sciences}.
\newblock Springer-Verlag, Berlin, 2005.
\newblock Invariant Theory and Algebraic Transformation Groups, V.

\bibitem{Rittatore1998}
Alvaro Rittatore.
\newblock Algebraic monoids and group embeddings.
\newblock {\em Transform. Groups}, 3(4):375--396, 1998.

\bibitem{Rota1964}
Gian-Carlo Rota.
\newblock On the foundations of combinatorial theory. {I}. {T}heory of
  {M}\"{o}bius functions.
\newblock {\em Z. Wahrscheinlichkeitstheorie und Verw. Gebiete}, 2:340--368
  (1964), 1964.

\bibitem{Spiegel2001}
Eugene Spiegel.
\newblock On the automorphisms of incidence algebras.
\newblock {\em J. Algebra}, 239(2):615--623, 2001.

\bibitem{SO}
Eugene Spiegel and Christopher~J. O'Donnell.
\newblock {\em Incidence algebras}, volume 206 of {\em Monographs and Textbooks
  in Pure and Applied Mathematics}.
\newblock Marcel Dekker, Inc., New York, 1997.

\bibitem{Stanley1970}
Richard~P. Stanley.
\newblock Structure of incidence algebras and their automorphism groups.
\newblock {\em Bull. Amer. Math. Soc.}, 76:1236--1239, 1970.

\bibitem{Timashev1994}
Dmitry~A. Timashev.
\newblock A generalization of the {B}ruhat decomposition.
\newblock {\em Izv. Ross. Akad. Nauk Ser. Mat.}, 58(5):110--123, 1994.

\bibitem{Timashev}
Dmitry~A. Timashev.
\newblock {\em Homogeneous spaces and equivariant embeddings}, volume 138 of
  {\em Encyclopaedia of Mathematical Sciences}.
\newblock Springer, Heidelberg, 2011.
\newblock Invariant Theory and Algebraic Transformation Groups, 8.

\bibitem{Vinberg1995}
\`Ernest~Borisovich Vinberg.
\newblock On reductive algebraic semigroups.
\newblock In {\em Lie groups and {L}ie algebras: {E}. {B}. {D}ynkin's
  {S}eminar}, volume 169 of {\em Amer. Math. Soc. Transl. Ser. 2}, pages
  145--182. Amer. Math. Soc., Providence, RI, 1995.

\end{thebibliography}

\end{document}